
\input amstex.tex
\loadmsam
\loadmsbm
\loadbold
\input amssym.tex
\input amstex
\baselineskip=13pt plus 2pt
\documentstyle{amsppt}
\pagewidth{33pc}
\hsize=12.3cm

\magnification=1200
\overfullrule=0pt

\def\enddemos{\hfill$\square$\enddemo}
\def\hb{\hfil\break}
\def\n{\noindent}

\def\pii{\pmb{\pi}}

\def\F{\Bbb F}

\def\Z{\Bbb Z}

\def\1{\bold 1}
\def\v{\bold v}
\def\mx{\operatorname{max}}
\def\Vol{\operatorname{Vol}}

\topmatter
\title On the symmetric square. Unstable Twisted Characters\endtitle
\author Yuval Z. Flicker and Dmitrii Zinoviev \endauthor
\footnote"~"{\n Department of Mathematics, The Ohio State University,
231 W. 18th Ave., Columbus, OH 43210-1174;\hb
email: flicker\@math.ohio-state.edu.\hb
\indent INRIA, Projet PRISME, 2004, rt. des Lucioles -BP.93,
06902, Sophia Antipolis Cedex, France, and\hb
\indent Inst. Prob. Info. Transmission, Russian Academy of Sciences,
Bolshoi Karetnyi per. 19, GSP-4, Moscow, 101447, Russia;\qquad
email: dimitrii.zinoviev\@sophia.inria.fr.\hb
\indent 1991 Mathematics Subject Classification: 10D40, 10D30, 12A67,
12A85, 14G10, 22E55, 11F27, 11R42, 11S40.}

\abstract We provide a purely local computation of the
(elliptic) twisted (by ``transpose-inverse'') character of the
representation $\pi=I(\1)$ of PGL(3) over a $p$-adic
field induced from the trivial representation of the maximal
parabolic subgroup. This computation is independent of the theory
of the symmetric square lifting of [IV] of automorphic and admissible
representations of SL(2) to PGL(3). It leads -- see [FK] --
to a proof of the (unstable) fundamental lemma in the theory of
the symmetric square lifting, namely that corresponding spherical
functions (on PGL(2) and PGL(3)) are matching: they have matching
orbital integrals. The new case in [FK] is the unstable one.
A direct local proof of the fundamental lemma is given in [V].
\endabstract
\endtopmatter

\document


This work continues the paper [FK], whose notations we use.
Our aim is to prove Proposition 1 of [FK] without using Theorem 0
there. Namely we provide a purely local computation of the twisted
character of $\pi=I(\1)$. Our model of $\pi$ is that of
[FK], where the twisted character $\chi_\pi$ is computed directly
and locally but only for the anisotropic twisted conjugacy class
$\delta'$ (see [FK], proof of Proposition 1). In [FK] the value
on the isotropic twisted conjugacy class $\delta$ is deduced from
the global Proposition 2.4 of [IV] -- recorded in [FK] as Theorem 0
-- which asserts that $\chi_\pi(\delta)=-\chi_\pi(\delta')$.

While the proof of Proposition 2.4 of [IV] is independent of the
results of [FK] (Theorems 1, 2, 3, 3$'$, which follow from Proposition
1), it is global, and so might lead some readers to worry that a
vicious circle is created. Moreover, the proof of this global result
requires heavy machinery. Here we provide a purely local proof of
Proposition 1 of [FK], and consequently make the results of [FK]
independent of Proposition 2.4 of [IV] (= Theorem 0 of [FK]).

Of course the conventional approach is to deduce the character
computation of [FK], Proposition 1, on using the global trace
formula comparison ([IV]) which is based on the fundamental lemma,
proven purely locally in [V]. The novel approach of [FK] -- which
we complete here -- is in reversing this perspective, and using
the global trace formula to prove the (unstable) fundamental lemma
from a purely local computation of the twisted character in a
special case.

Further, an independent, direct computation of the very precise
character calculation gives another assurance of the validity of
the trace formula approach to the lifting project. It will be
interesting to develop this approach in other lifting situations,
especially since our technique is different from the well-known,
standard techniques of trace formulae and dual reductive pairs.
A first step in this direction was taken in our work [FZ], where
the twisted -- by the transpose-inverse involution -- character
of a representation of PGL(4) analogous to the one considered here,
is computed. The situation of [FZ] is new, dealing with the
exterior product of two representations of GL(2) and the structure
of representations of the rank two symplectic group. Such character
computations are not yet available by any other technique. However,
the computations of [FZ] -- although elementary -- are involved,
as they depend on the classification of [F] of the twisted conjugacy
classes in GL(4). This is another reason for the present work, which
considers the initial non trivial case of our technique -- where the
computations are still simple and can clarify the method. We believe
that our methods, pursued in [FZ] in a more complicated case, would
apply in quite general lifting situations, in conjunction with, and
as an alternative to the trace formula.

\bigskip

Proposition 1 of [FK] asserts that if $\1$ is the trivial
PGL(2,$F$)-module, $\pi=I(\1)$ is the PGL(3,$F$)-module
normalizedly induced from the trivial representation of the maximal
parabolic subgroup (whose Levi component is GL(2,$F$)), and $\delta$
is a $\sigma$-regular element of PGL(3,$F$) with elliptic regular
norm $\gamma_1=N_1\delta$, then $(\Delta(\delta)/\Delta_1(\gamma_1))
\chi_\pi(\delta)=\kappa(\delta)$.

The proof of Proposition 1 in [FK] reduces this to the claim that the
value at $s=0$ of
$$|4u\theta|^{1/2}|u(\alpha^2-\theta)|^{-s/2}
\int_{V^0}\,|x^2+uy^2-\theta z^2|^{3(s-1)/2}\,dx\,dy\,dz$$
is $-\kappa(\delta)q^{-1/2}(1+q^{-1/2}+q^{-1})$ (see bottom of page
499, and Lemma 2, in [FK]).

This equality is verified in [FK], p. 499,
when the quadratic form $x^2+uy^2-\theta z^2$ is anisotropic, in
which case $\kappa(\delta)=-1$ and the integral converges for all $s$.

Here we deal with the case where the quadratic form is {\bf isotropic},
in which case $\kappa(\delta)=1$, the integral converges only in some
half plane of $s$, and the value at $s=0$ is obtained by analytic
continuation.

Recall that $F$ is a local non-archimedean field of odd residual
characteristic; $R$ denotes the (local) ring of integers of $F$;
$\pii$ signifies a generator of the maximal ideal of $R$. Denote
by $q$ the number of elements of the residue field $R/\pii R$ of $R$.
By  $\F$ we mean a set of representatives in $R$ for the
finite field $R/\pii$. The absolute value on $F$ is normalized by
$|\pii|=q^{-1}$.

The case of interest is that where $K=F(\sqrt\theta)$ is a quadratic
extension of $F$, thus $\theta\in F^\times-F^{\times 2}$. Since the
twisted character depends only on the twisted conjugacy class, we may
assume that $|\theta|$ and $|u|$
lie in $\{1,q^{-1}\}$.

\proclaim{0. Lemma} We may assume that the quadratic form
$x^2 + u y^2 -\theta z^2$ takes one of three avatars: $x^2-\theta z^2-y^2$,
$\theta\in R-R^2$; $x^2-\pii z^2+\pii y^2$; or $x^2-\pii z^2-y^2$.
\endproclaim

\demo{Proof}
(1) If $K/F$ is unramified, then $|\theta|=1$, thus
$\theta\in R^\times-R^{\times 2}$. The norm group
$N_{K/F}K^\times$ is $\pi^{2\Z}R^\times$. If
$x^2-\theta z^2+uy^2$ represents 0 then $-u\in R^\times$.
If $-1$ is not a square, thus $\theta=-1$, then $u$ is
$-1$ (get $x^2-z^2-y^2$) or $u=1$ (get $x^2-z^2+y^2$,
equivalent case). If $-1\in R^{\times 2}$, the case
of $u=\theta$ ($x^2-\theta z^2+\theta y^2=\theta(y^2
+\theta^{-1}x^2-z^2$)) is equivalent to the case of
$u=-1$. So wlog $u=-1$ and the form is $x^2-\theta z^2-y^2$,
$|u\theta|=1$.

\n (2) If $K/F$ is ramified, $|\theta|=q^{-1}$ and
$N_{K/F}K^\times=(-\theta)^{\Z}R^{\times 2}$. The form
$x^2-\theta z^2+uy^2$ represents zero when $-u\in
R^{\times 2}$ or $-u\in -\theta R^{\times 2}$. Then the
form looks like $x^2-\theta z^2+\theta y^2$ with $u=\theta$ and
$|\theta u|=q^{-2}$, or $x^2-\theta z^2-y^2$ with $u=-1$ and
$|\theta u|=q^{-1}$.
The Lemma follows.
\enddemos

\medskip
We are interested in the value at $s=-3/2$ of the integral
$I_s(u,\theta)$ of $|x^2 + u y^2 - \theta z^2|^s$ over the set
$V^0=V/\sim$, where
$$V=\{\v=(x,y,z)\in R^3;\mx\{|x|,|y|,|z|\}=1\}$$
and $\sim$ is the equivalence relation $\v\sim\alpha \v$ for
$\alpha\in R^\times$.

The set $V^0$ is the disjoint union of the subsets
$$V_n^0=V_n^0(u,\theta)=V_n(u,\theta)/\sim,$$ where
$$V_n=V_n(u,\theta) = \{\v; \mx\{|x|,|y|,|z|\}=1,
|x^2 + u y^2 -\theta z^2|=1/q^n\},$$
over $n\ge 0$, and of $\{\v; x^2+uy^2-\theta z^2=0\}/\sim$,
a set of measure zero.

Thus we have
$$I_s(u,\theta) = \sum_{n=0}^{\infty}q^{-ns}\Vol(V_n^0(u,\theta)).$$
\medskip

\proclaim{Theorem} The value of $|u\theta|^{1/2}I_s(u,\theta)$
at $s=-3/2$ is $-q^{-1/2}(1+q^{-1/2}+q^{-1})$.
\endproclaim

The problem is simply to compute the volumes
$$\Vol(V_n^0(u,\theta))=\Vol(V_n(u,\theta))/(1-1/q)\qquad (n\ge 0).$$

\proclaim{1. Lemma} When $\theta=\pii$ and $u=-1$, thus $|u\theta|=1/q$,
we have
$$\Vol(V_n^0)= \cases
(1-1/q),              & \qquad\text{if}\ n=0, \\
2q^{-1}(1-1/q)+1/q^2, & \qquad\text{if}\ n=1, \\
2q^{-n}(1-1/q),       & \qquad\text{if}\ n\ge 2.
\endcases$$
\endproclaim

\demo{Proof} Recall that
$$V_0=V_0(-1,\pii)=\{(x,y,z);\,\mx\{|x|,|y|,|z|\}=1,\,
|x^2-y^2-\pii z^2|=1\}.$$
Since $|z|\le 1$, we have $|\pii z^2| < 1$, and
$$1 = |x^2 - y^2 -\pii z^2|=|x^2 - y^2|=|x-y||x+y|.$$
Thus $|x-y|=|x+y|=1$, and if $|x|\neq |y|$, $|x\pm y|=\mx\{|x|,|y|\}$.
We split $V_0$ into three distinct subsets, corresponding to the cases
$|x|=|y|=1$; $|x|=1$, $|y|<1$; and $|x|<1$, $|y|=1$.
The volume is then
$$\Vol(V_0)=\int_{|z|\le
1}\int_{|x|=1}\left[\int_{|y|=1,|x-y|=|x+y|=1}
\right]dydxdz$$
$$ + \int_{|z|\le
1}\left[\int_{|x|=1}\int_{|y|<1}+\int_{|x|<1}\int_{|y|=1}
\right]dydxdz$$
$$ = \int_{|x|=1}\left[\int_{|y|=1,|x-y|=|x+y|=1}\right]dydx
+ \frac{2}{q}\left(1-\frac{1}{q}\right)
=\left(1-\frac{1}{q}\right)^2.$$

To consider the $V_n$ with $n\ge 1$, where
$|x^2 - y^2 -\pii z^2| = 1/q^n$, recall that any $p$-adic number $a$
such that $|a|\le 1$ can be written as a power series in $\pii$:
$$a = \sum_{i=0}^{\infty}a_i\pii^i = a_0+a_1\pii+a_2\pii^2+...\qquad
(a_i\in\F).$$
In particular $|a|=1/q^n$ implies that $a_0=a_1=...=a_{n-1}=0$,
and $a_n\ne 0$. If
$$x = \sum_{i=0}^{\infty}x_i\pii^i,\qquad
y = \sum_{i=0}^{\infty}y_i\pii^i,\qquad
z = \sum_{i=0}^{\infty}z_i\pii^i\qquad (x_i,y_i,z_i\in\F),$$
then
$$x^2 = \sum_{i=0}^{\infty}a_i\pii^i,\qquad
y^2 = \sum_{i=0}^{\infty}b_i\pii^i,\qquad
z^2 = \sum_{i=0}^{\infty}c_i\pii^i,$$
where
$$
a_i = \sum_{j=0}^{i}x_j x_{i-j},\qquad b_i = \sum_{j=0}^{i}y_j y_{i-j},\qquad
c_i = \sum_{j=0}^{i}z_j z_{i-j}\qquad (a_i,b_i,c_i\in\F).$$

We have
$$x^2 - y^2 -\pii z^2 = \sum_{i=0}^{\infty}f_i\pii^i\qquad
(f_i\in\F),$$
where $f_0=a_0-b_0,\ f_i=a_i - b_i - c_{i-1}\ (i\ge 1)$.
Since $|x^2 - y^2 -\pii z^2| = 1/q^n$, we have that
$f_0=f_1=...=f_{n-1}=0$, and $f_n\ne 0$.
Thus we obtain the relations (for $a$, $b$, $c$ in the set $\F$,
which (modulo $\pii$) is the field $R/\pii$):
$$a_0-b_0=0,\qquad a_i - b_i - c_{i-1}=0\qquad (i=1,...,n-1),\qquad
a_n - b_n - c_{n-1}\ne 0.$$
Recall that together with $\mx\{|x|,|y|,|z|\}=1$, these relations define
the set $V_n$.

To compute the volume of $V_n$ we integrate in the order: $...dydzdx$.
From $a_0 - b_0=0$ it follows that $y_0=\pm x_0$, and from
$a_i - b_i - c_{i-1}$ ($i\ge 1$) it follows that
$$2y_0y_i = a_i- c_{i-1} - \sum_{j=1}^{i-1}y_j y_{i-j},$$
where in the case of $i=1$ the sum over $j$ is empty.

Let $n\ge 2$. When $i=1$ we have $2x_0x_1-2y_0y_1-z_0^2=0$.
So if $x_0=0$ (in $R/\pii$, i.e. $|x|<1$), it follows that $y_0=0$ and
$z_0=0$ (i.e. $|y|<1$, $|z|<1$). This contradicts the fact
that $\mx\{|x|,|y|,|z|\}=1$. Thus $|x|=1$. In this case $y_0\ne 0$
and (for $n\ge 2$) we have:
$$\Vol(V_n)=\int_{|x|= 1}\int_{|z|\le 1}\left[\int dy\right]dzdx,$$
where the variable $y$ is such that once written as
$y=y_0+y_1\pii+y_2\pii^2+...$, it has to satisfy: $y_0=\pm x_0$, and
$y_i$ ($i=1,...,n-1$) is defined uniquely from $a_i-b_i-c_{i-1}=0$,
and $y_n\ne$ some value defined by $a_n - b_n - c_{n-1}\ne 0$.
Thus when $n\ge 2$,
$$\Vol(V_n) = \frac{2}{q}\left(\frac{1}{q}\right)^{n-1}
\left(1-\frac{1}{q}\right)^2 = \frac{2}{q^n}\left(1-\frac{1}{q}\right)^2.$$

Let $n=1$. When $i=1$ we have $2x_0x_1-2y_0y_1-z_0^2\ne 0$.
So if $x_0=0$ (i.e. $|x|<1$), it follows that $y_0=0$ and
$z_0\ne 0$ (i.e. we have an additional contribution from
$|x|<1$, $|y|<1$, $|z|=1$). Thus,
$$\Vol(V_1) = \frac{2}{q}\left(1-\frac{1}{q}\right)^2 +
\frac{1}{q^2}\left(1-\frac{1}{q}\right).$$
The Lemma follows.
\enddemos

\proclaim{2. Lemma} When $u$ and $\theta$ equal $\pii$, thus
$|u\theta|=1/q^2$, we have
$$\Vol(V_n^0)=\cases
1,               & \qquad\text{if}\ n=0, \\
q^{-1}(1-1/q),   & \qquad\text{if}\ n=1, \\
2q^{-n}(1-1/q),  & \qquad\text{if}\ n\ge 2.
\endcases $$
\endproclaim

\demo{Proof} To compute $\Vol(V_0)$, recall that
$$V_0 = \{(x,y,z); \mx\{|x|,|y|,|z|\}=1, |x^2 +\pii (y^2 - z^2)|=1\}.$$
Since $|y|\le 1$, $|z|\le 1$, we have $|x^2 +\pii (y^2 - z^2)|=|x^2|=1$,
and so
$$\Vol(V_0)=\int_{|z|\le 1}\int_{|y|\le 1}\int_{|x|=1}dxdydz
=1-\frac{1}{q}.$$

To compute $\Vol(V_n)$, $n\ge 1$, recall that
$$V_n = \{(x,y,z); \mx\{|x|,|y|,|z|\}=1,
|x^2 +\pii (y^2 - z^2)|=1/q^n\}.$$
Following the notations of Lemma 1 we write
$$x^2 + \pii(y^2  -z^2) = \sum_{i=0}^{\infty}f_i\pii^i\qquad
(f_i\in\F),$$
where $f_0=a_0$ and $f_i=a_i + b_{i-1} - c_{i-1}$ ($i\ge 1$).
The condition which defines $V_n$ is that $f_0=f_1=...=f_{n-1}=0$
and $f_n\ne 0$. The equation $f_0=0$ implies that $x_0=0$ (i.e. $|x|<1$).
We arrange the order of integration to be: $...dydzdx$.

When $n\ge 2$, since $x_0=0$, $f_1=0$ implies
that $y_0^2-z_0^2= 0$. Using $\mx\{|x|,|y|,|z|\}=1$ we conclude that
$y_0=\pm z_0\ne 0$ (i.e. $|z|=1$, $|z^2-y^2|<1$). Thus we have
$$\Vol(V_n)=\int_{|x|<1}\int_{|z|= 1}\left[\int dy\right]dzdx$$
where the variable $y$ is such that once written as
$y=y_0+y_1\pii+y_2\pii^2+...$, it has to satisfy: $y_0=\pm z_0$, and
$y_i$ ($i=1,...,n-2$) is defined uniquely from $a_i+b_{i-1}-c_{i-1}=0$,
and $y_{n-1}\ne$ some value defined by $a_n + b_{n-1} - c_{n-1}\ne 0$.
Thus when $n\ge 2$,
$$\Vol(V_n) = \frac{1}{q}\frac{2}{q}\left(\frac{1}{q}\right)^{n-2}
\left(1-\frac{1}{q}\right)^2 = \frac{2}{q^{n}}\left(1-\frac{1}{q}\right)^2.$$

When $n=1$ we have $f_0=0$, $f_1\ne 0$. These amount to $x_0=0$,
$y_0\ne\pm z_0$. Separating the two cases $z_0=0$, and
$z_0\ne 0$, we obtain
$$\Vol(V_1)=\int_{|x|<1}\int_{|z|<1}\int_{|y|=1}dydzdx +
\int_{|x|<1}\int_{|z|=1}\int_{|y^2-z^2|=1}dydzdx$$
$$=\frac{1}{q^2}\left(1-\frac{1}{q}\right)+
\frac{1}{q}\left(1-\frac{1}{q}\right)\left(1-\frac{2}{q}\right)=
\frac{1}{q}\left(1-\frac{1}{q}\right)^2.$$
The Lemma follows.
\enddemos
\medskip

\proclaim{3. Lemma} When $K/F$ is unramified, thus
$|u\theta|=1$, we have
$$\Vol(V_n^0)=\cases
1,                    & \qquad\text{if}\ n=0, \\
q^{-n}(1-1/q)(1+1/q), & \qquad\text{if}\ n\ge 1.
\endcases$$
\endproclaim

\demo{Proof} First we compute $\Vol(V_0)$. Recall that
$$V_0 = \{(x,y,z); \mx\{|x|,|y|,|z|\}=1, |x^2 - y^2 - \theta z^2|=1\}.$$
Since $|x^2 - y^2 - \theta z^2|\le \mx\{|x|,|y|,|z|\}$,
$$V_0 = \{(x,y,z)\in R^3; |x^2 - y^2 - \theta z^2|=1\}.$$
Making the change of variables $u=x+y$, $v=x-y$, we obtain
$$V_0 = \{(u,v,z)\in R^3; |uv - \theta z^2|=1\}.$$

Assume that $|uv|<1$. Since $|uv - \theta z^2|=1$, it follows that
$|z|=1$. The contribution from the set $|uv|<1$ is
$$\int_{|z|=1}\left[\int_{|u|<1}\int_{|v|\le1}+\int_{|u|=1}\int_{|v|<1}
\right]dudvdz$$
$$ = \left(1-\frac{1}{q}\right)\left(\frac{1}{q}+
\left(1-\frac{1}{q}\right)\frac{1}{q}\right)
=\frac{1}{q}\left(1-\frac{1}{q}\right)\left(2-\frac{1}{q}\right).$$

Assume that $|uv|=1$, i.e. $|u|=|v|=1$. We arrange the order of
integration as: $dudvdz$. If $|z|<1$ then $|uv-\theta z^2|=|uv|=1$.
If $|z|=1$ we introduce $U(v,z)= \{u;\ |u|=1,\ |uv-\theta z^2|=1\}$,
a set of volume $1-2/q$, and note that the contribution from
the set $|uv|=1$ is
$$\int_{|z|<1}\int_{|v|=1}\int_{|u|=1}dudvdz +
\int_{|z|=1}\int_{|v|=1}\int_{U(v,z)}dudvdz.$$
The sum of the two integrals is
$$\frac{1}{q}\left(1-\frac{1}{q}\right)^2+\left(1-\frac{1}{q}\right)^2
\left(1-\frac{2}{q}\right)=\left(1-\frac{1}{q}\right)^3.$$

Adding the contributions from $|uv|<1$ and $|uv|=1$ we then obtain
$$\Vol(V_0) = \frac{1}{q}\left(1-\frac{1}{q}\right)
\left(2-\frac{1}{q}\right) + \left(1-\frac{1}{q}\right)^3
= 1-\frac{1}{q}.$$

Next we compute $\Vol(V_n)$, $n\ge 1$. Recall that
$$V_n = \{(x,y,z); \mx\{|x|,|y|,|z|\}=1,
|x^2 - y^2 - \theta z^2|=1/q^n\}.$$
Making the change of variables $u=x+y$, $v=x-y$, we obtain
$$V_n = \{(u,v,z); \mx\{|u+v|,|u-v|,|z|\}=1,|uv - \theta z^2|=1/q^n\}.$$
Since the set $\{v=0\}$ is of measure zero, we assume that
$v\ne 0$. Then $|uv - \theta z^2|=1/q^n$ implies that
$u=\theta z^2 v^{-1} + tv^{-1}\pii^n$, where $|t|=1$.
There are two cases.

Assume that $|v|=1$. Note that if
$|z|=1$, then $\mx\{|u+v|,|u-v|,|z|\}=1$ is satisfied, and if
$|z|<1$, then (recall that $n\ge 1$)
$$|u|=|\theta z^2 v^{-1} + tv^{-1}\pii^n|\le\mx\{|z^2|,q^{-n}\}<1,$$
and $|u+v|=|v|=1$. So $|v|=1$ implies $\mx\{|u+v|,|u-v|,|z|\}=1$. Further,
since $|v|=1$, we have $du=q^{-n}dt$. Thus the contribution from the set
with $|v|=1$ is
$$\int_{|z|\le 1}\int_{|v|=1}\int_{|uv-\theta z^2|=1/q^n}dudvdz
=\int_{|z|\le 1}\int_{|v|=1}\int_{|t|=1}\frac{dt}{q^n}dvdz
=\frac{1}{q^n}\left(1-\frac{1}{q}\right)^2.$$

Assume that $|v|<1$. Note that if $|z|=1$, since $|u|\le 1$ we have
$q^{-n}=|uv-\theta z^2|=|\theta z^2|=1$, a contradiction.
Thus $|z|<1$, and in order to satisfy $\mx\{|u+v|,|u-v|,|z|\}=1$,
we should have $|u|=1$. The contribution from the set with $|v|<1$ is
$$\int_{|z|<1}\int_{|u|=1}\int_{|uv-\theta z^2|=1/q^n}dvdudz.$$
We write $v=\theta z^2 u^{-1} + tu^{-1}\pii^n$, where $|t|=1$,
and $dv=q^{-n}dt$. The integral equals
$$\int_{|z|<1}\int_{|u|=1}\int_{|t|=1}\frac{dt}{q^n}dudz
=\frac{1}{q}\frac{1}{q^n}\left(1-\frac{1}{q}\right)^2.$$

Adding the contributions from $|v|=1$ and $|v|<1$ we obtain
$$\Vol(V_n)=\frac{1}{q^n}\left(1-\frac{1}{q}\right)^2+
\frac{1}{q}\frac{1}{q^n}\left(1-\frac{1}{q}\right)^2
= \frac{1}{q^n}\left(1-\frac{1}{q}\right)^2\left(1+\frac{1}{q}\right).$$
The Lemma follows.
\enddemos

This completes the proof of the theorem, so that we provided a purely
local proof of (the character relation of) Proposition 1 of [FK]. We
believe that analogous computations can be carried out in other lifting
situations, to provide direct and local computations of twisted
characters. As noted in the introduction, a step in this direction
is taken in [FZ].

\bigskip

\def\refe#1#2{\n\hangindent 5em\hangafter1\hbox to 5em{\hfil#1\quad}#2}
\subheading{References}
\medskip

\refe{[IV]}{Y. Flicker, On the symmetric square. Applications of a trace
formula, {\it Trans. AMS} 330 (1992), 125-152.}

\refe{[V]}{Y. Flicker, On the symmetric-square. The fundamental lemma, {\it
Pacific J. Math.} 175 (1996), 507-526.}

\refe{[F]}{Y. Flicker, Matching of orbital integrals on GL(4) and GSp(2),
{\it Mem. Amer. Math. Soc.}  137 (1999), 1-114.}

\refe{[FK]}{Y. Flicker, D. Kazhdan, On the symmetric-square: Unstable
local transfer, {\it Invent. Math.} 91 (1988), 493-504.}

\refe{[FZ]}{Y. Flicker, D. Zinoviev, Twisted character of a small
representation of PGL(4), Moscow Mathematical Journal, V. 4,
N$^{\circ}$ 2, April-June 2004, 333--368. }
\enddocument